\documentclass[11pt]{article}

\usepackage{amsmath,amsthm,amsfonts}

\usepackage{amssymb}

\usepackage{epsfig}

\usepackage{rotating}

\usepackage{subfigure}

\newcommand{\p}[2]{\frac{\partial#1}{\partial#2}}
\newcommand{\pd}[2]{\frac{\partial#1}{\partial#2}}
\begin{document}

\title{Renormalized reduced models for singular PDEs}
\author{Panos Stinis \\ 
Department of Mathematics \\
University of Minnesota \\
    Minneapolis, MN 55455} 

\date {}

\maketitle

\begin{abstract}
We present a novel way of constructing reduced models for systems of ordinary differential equations. In particular, the approach combines the concepts of renormalization and effective field theory developed in the context of high energy physics and the Mori-Zwanzig formalism of irreversible statistical mechanics. The reduced models we construct depend on coefficients which measure the importance of the different terms appearing in the model and need to be estimated. The proposed approach allows the estimation of these coefficients on the fly by enforcing the equality of integral quantities of the solution as computed from the original system and the reduced model.  In this way we are able to construct stable reduced models of higher order than was previously possible. The method is applied to the problem of computing reduced models for ordinary differential equation systems resulting from Fourier expansions of singular (or near-singular) time-dependent partial differential equations. Results for the 1D Burgers and the 3D incompressible Euler equations are used to illustrate the construction. We also present, for the 1D Burgers and the 3D Euler equations, a simple and efficient recursive algorithm for calculating the higher order terms.
\end{abstract}

%\maketitle

\section*{Introduction}

Spatial discretizations or Fourier expansions of the solutions of partial differential equations (PDEs) which depend on time lead to systems of ordinary differential equations (ODEs). The most difficult case arises when the solution of a PDE becomes singular in finite time. At such instants the solution of the PDE develops activity down to the zero length scale. A brute force numerical simulation (no matter how large) of such a solution is bound to fail because the simulation has a finite resolution and thus will be unable to resolve all the length scales down to the zero scale. When the solution develops activity at a scale smaller than the smallest scale available to the simulation, the numerically computed solution becomes underresolved. This leads to a rapid deterioration of the accuracy of the simulation.

The notion of propagation of activity to smaller and smaller scales depends on the physical context of the PDE. In some cases, e.g. the 3D Euler or Navier-Stokes equations \cite{doering}, this could mean a cascade of energy to smaller and smaller scales. In other cases, e.g. the nonlinear focusing Schr\"odinger equation \cite{sulem}, this could mean a cascade of mass to smaller and smaller scales. Irrespective of the specific physical context, the problem facing the numerical analyst is how to use a {\it finite} simulation and yet prevent the computed solution from suffering a loss of accuracy. In other words, how to construct a numerical method which reproduces correctly the features of the solution of the original equation at  the length scales that are available numerically. This is the motivation behind the construction of reduced models (see e.g. \cite{givon,CS05}).   

By construction, a reduced model must allow for energy (mass) to escape from the scales that are accessible to the simulation (called resolved scales or modes) to the inaccessible scales (called unresolved). The main difficulty in constructing an accurate reduced model is the need to estimate the {\it correct} rate at which activity is propagated from the resolved to the unresolved scales. The Mori-Zwanzig (MZ) formalism \cite{CHK00,CHK3}  proceeds by dividing the available resolution into resolved and unresolved parts. Then, it constructs a reduced model for the resolved scales and uses the unresolved scales to effect the drain of energy (mass) out of the resolved scales.  

Although the MZ formalism allows for the construction, in principle, of an exact reduced model it has two drawbacks (which are also shared by {\it any} other reduction formalism). First, the reduced model can be, in general, prohibitively expensive to calculate. The reason is that one must obtain an accurate representation of the behavior of the unresolved scales before they can be safely eliminated. Obtaining this representation can be rather costly.

The second drawback is more subtle and has not been adequately appreciated by the scientific computing community. It is specific to the case of constructing reduced models for singular PDEs or in general for systems of ordinary differential equations which are larger than any available numerical resolution. Suppose that you have to construct a reduced model of a full system which is larger than any available numerical simulation. Let us call this system S1. Exactly because S1 is larger than any available numerical simulation, if we want to construct a reduced model we have to use as a starting point a system, call it S2, whose size is smaller than the size of S1. Suppose that you start with S2 and use the MZ formalism (or any other reduction formalism for that matter) and construct an {\it exact} reduced model S3 for a subset of S2. An exact reduced model means that if one evolves S2 and S3 separately, then the behavior of the scales resolved by S3 will be the same as the behavior of the scales in S3 predicted by the simulation of the system S2. However, and this is the heart of the problem, since S2 itself will become eventually underresolved, the exact reduced model S3 will also become underresolved. In other words, the predictions of the exact reduced model S3 can only be trusted for as long as the predictions of the system S2 can be trusted. As a result, {\it any reduced model that has any chance of being accurate for longer times cannot be exact.}

There are examples of {\it inexact} reduced models, e.g. the $t$-model \cite{CHK3,bernstein,HS06,S10}, coming from the MZ formalism, which have been applied to singular PDEs and shown numerically to be relatively accurate for long time intervals. However, the $t$-model's accuracy is difficult to assess beforehand and the reason for its relative success has remained a mystery (see also \cite{chandy} for an application to the 3D Navier-Stokes equations which shows that the $t$-model, while not bad, is in need of some modification). In order to construct better reduced models we need to incorporate dynamic information from the full system which will help us decide which of the terms appearing in the exact reduced model are the ones that are most important. In this way, we can construct an inexact but accurate reduced model by keeping the important terms and disregarding the unimportant ones. The underlying principle of this approach is the same as the principle behind the construction of effective field theories in high energy physics \cite{georgi}. 

The way we propose to address the problem of constructing better reduced models is to embed the MZ reduced models in a larger class of reduced models which share the same functional form as the MZ reduced models but have different coefficients in front of the various terms that appear in the reduced models. Then, one can estimate these coefficients on the fly while the original system of equations is still valid. The estimation of the coefficients is achieved by requiring that certain integral quantities (e.g. $l_p$ norms) involving only resolved scales, should acquire the same values when computed from the original system and the reduced model. The constraints used to obtain the coefficients are the analog of the ``matching conditions" used in effective field theory. Also, the approach is the time-dependent analog of the process of renormalization used in high energy and condensed matter physics \cite{collins,goldenfeld}.  Before the original system ceases to be valid, one reverts to the reduced model with the various coefficients having their estimated values. We call the proposed approach the renormalized Mori-Zwanzig algorithm.  

In addition to computing the coefficients of the reduced model on the fly, the proposed approach allows us to gauge the importance of the different terms appearing in the reduced model. While at the moment a rigorous proof is lacking, the numerical results (see Section \ref{universality}) suggest that the values of the coefficients depend on the ratio of the smallest scale present in the initial condition to the smallest scale resolved by the reduced model. The larger is this ratio, the larger are the coefficients for more terms and thus, more terms are needed in the reduced model. On the other hand, for the case when the ratio is small only a few number of terms need to be kept without sacrificing the accuracy of the model. As we have already said, this situation is analogous to the principle underlying the effective field theory method in high energy physics \cite{georgi}. It is also interesting to investigate to what extent the values of the coefficients are determined by the scaling symmetries of the underlying PDE. This could allow us to identify different universality classes of PDEs based on their scaling symmetries (see also Section \ref{universality}).

Finally, we note that the proposed method of computing the coefficients was first presented by the author in \cite{S09}. The goal in \cite{S09} was to construct a mesh refinement scheme to allow us to reach the singularity instant more efficiently. Then, the computed coefficients were used in a fixed point analysis of the singularity in order to estimate the blow-up exponents. In the current work the goal is different. We not only want to reach the singularity instant or, more generally, the time that the full system runs out of resolution, but also follow the solution for later times. Also, in the current work, we present a new way of expanding the memory term of the MZ formalism. This allows us, at least for Burgers and Euler, to calculate recursively and efficiently (and with minimal storage requirements) the higher order terms in the memory expansion. This is a welcome feature since analytical computation of the higher order terms quickly becomes exhausting and, thus,  prone to errors.

The paper is organized as follows. Section \ref{algo} presents the main construction along with a brief presentation of the MZ formalism and how to compute the coefficients of the reduced model on the fly. Section \ref{numerics} presents numerical results for the 1D Burgers and the 3D incompressible Euler equations.  Section \ref{conclusions} contains a discussion of the results and presents some directions for future work.

%%%%%%%%%%%%%%%%End of Introduction%%%%%%%%%%%%%%%%%%%%%

\section{Renormalization of Mori-Zwanzig reduced models}\label{algo}

In this section we explain the main ideas behind the proposed approach. In Section \ref{full_reduced} we set up the notation for the original system and the reduced model in an abstract way which does not make reference to any specific method for obtaining the reduced model. In Section \ref{coefficient_compute} we show how to obtain the coefficients for the reduced model. In Section \ref{mzformalism} we give a brief presentation of the MZ formalism which allows us to obtain the functional form of the terms appearing in the reduced model. In Section \ref{combine} we combine the ideas in Section \ref{coefficient_compute} with the MZ formalism from Section \ref{mzformalism} to derive the proposed algorithm for computing renormalized MZ reduced models.

\subsection{Full and reduced systems}\label{full_reduced}

Suppose that we want to construct a reduced model for the partial differential equation (PDE)
$$ v_t + H (t,x,v,v_x,...)=0 $$
where $H$ is a, in general nonlinear, operator and $x \in D \subseteq \mathbb{R}^d$ (the construction extends readily to the case of systems of partial differential equations). After spatial discretization or expansion of the solution in series, the PDE transforms into a system of ordinary differential equations (ODEs). For simplicity we restrict ourselves to the case of periodic boundary conditions, so that a Fourier expansion of the solution leads to system of ODEs for the Fourier coefficients. To simulate the system for the Fourier coefficients we need to truncate at some point the Fourier 
expansion. Let $F \cup G$ denote the set of Fourier modes retained in the series, where we have split the Fourier modes in two sets, $F$ and $G.$ 
We call the modes in $F$ resolved and the modes in $G$ unresolved. The reduced model involving only the resolved modes $F$ will be called the reduced system and the system involving both the resolved {\it and} unresolved modes $F \cup G$ will be called the full system. 

The main idea behind the algorithm is that the evolution of moments of the reduced set of modes, for example $l_p$ norms of the modes in $F$, should be the same whether computed from the full or the reduced system. This requirement will eventually allow us to compute the coefficients appearing in the reduced model (see Section \ref{coefficient_compute}). 

The full system of equations for the modes $F \cup G$ is given by 
$$\frac{du(t)}{dt} = R (t,u(t)),$$
where $u = ( \{u_k\}), \; k \in F \cup G$ is the vector of Fourier coefficients of $u$ and $R$ is the Fourier transform of the operator $H.$ The system should be supplemented with an initial condition $u(0)=u_0.$ The vector of Fourier coefficients can be written as $ u = (\hat{u}, 
\tilde{u}),$ where $ \hat{u}$ are the resolved modes (those in $F$) and $\tilde{u}$ the unresolved ones (those in $G$). Similarly, for the right hand sides (RHS) we have $R(t,u) = (\hat{R}(t,u), \tilde{R}(t,u)).$ Note that the RHS of the resolved modes involves both resolved and unresolved modes. In anticipation of the construction of a reduced model we can rewrite the RHS as $R(t,u)=R^{(0)}(t,u) = (\hat{R}^{(0)}(t,u), \tilde{R}^{(0)}(t,u)).$  

In general, when one constructs a reduced model, additional terms appear on the RHS of the equations of the reduced model (see Section \ref{mzformalism} for more details). The role of these additional terms is to account for the interactions between the resolved and unresolved modes, since the unresolved modes no longer appear explicitly in the reduced model. As is standard in renormalization theory \cite{binney}, one can augment the RHS of the equations in the full system by including such additional terms. That is accomplished by multiplying each of these additional terms by a zero coefficient. In this way, the reduced and full systems' RHSs have the same functional form. In particular, for each mode $u_k, \; k \in F \cup G,$ we can rewrite $R_k^{(0)}(t,u)$ as 
$$R_k^{(0)}(t,u(t)) = \sum_{i=1}^{m} a^{(0)}_i R^{(0)}_{ik}(t,u(t)),$$
where $R^{(0)}_{1k}(t,u(t))=R_k^{(0)}(t,u(t))$ and $R^{(0)}_{ik}(t,u(t)), \; \text{for} \;  i=2,\ldots,m$ are of the same functional form as the additional terms which appear in the reduced model. This is easy to do by taking $ a^{(0)}_1=1$ and $a^{(0)}_i=0, \; \text{for} \;  i=2,\ldots,m.$ Thus, the equation for the the mode $u_k, \; k \in F \cup G$ is written as 
\begin{equation}\label{full}
\frac{d{u_k}(t)}{dt} = {R}_k(t,u)={R}_k^{(0)} (t,u(t))=\sum_{i=1}^{m} a^{(0)}_i {R}^{(0)}_{ik} (t,u(t))
\end{equation}
Correspondingly, the reduced model for the mode $u'_k, \; k \in F $ is given by 
\begin{equation}\label{reduced}
\frac{d{u}_k'(t)}{dt} = R_k^{(1)} (t,\hat{u}'(t))=\sum_{i=1}^{m} a^{(1)}_i R^{(1)}_{ik} (t,\hat{u}'(t))
\end{equation}
with initial condition ${u}_k'(0)={u}_{0k}.$ We repeat that the functions $R^{(1)}_{ik}, \; i=1,\ldots,m, \; k \in F,$ have the same form as the functions $R^{(0)}_{ik}, \; i=1,\ldots,m, \; k \in F,$ but they depend {\it only} on the reduced set of modes $F.$ Dimensional reduction transforms the vector $a^{(0)}=(a^{(0)}_1,\ldots,a^{(0)}_m)$ to $a^{(1)}=(a^{(1)}_1,\ldots,a^{(1)}_m).$ This allows one to determine the relation of the full to the reduced system by focusing on the change of the vector $a^{(0)}$ to $a^{(1)}.$ Also, the vectors $a^{(0)}$ and $a^{(1)}$ do not have to be constant in time. This does not change the analysis that follows.

Define $m$ quantities 
$\hat{E}_i, \; i=1,\ldots,m$ involving only modes in $F.$ For example, these could be $l_p$ norms of the reduced set of modes. To proceed we require that these quantities' rates of change are the same when computed from (\ref{full}) and (\ref{reduced}), i.e. 
\begin{equation}\label{conditions}
\frac{d\hat{E}_i(\hat{u})}{dt} = \frac{d\hat{E}_i(\hat{u}')}{dt}, \; i=1,\ldots,m.
\end{equation}
Note that similar conditions, albeit time-independent, lie at the heart of the renormalization group theory for equilibrium systems (\cite{binney} p. 154). Also, the conditions \eqref{conditions} are the analog of the ``matching conditions" underlying the construction of effective field theories \cite{georgi}.

\subsection{How to compute the coefficients of the reduced model}\label{coefficient_compute}
When we only know the functional form of the terms appearing in the reduced model but not their coefficients it is not possible to evolve a reduced system. We present a way of actually computing the coefficients of the reduced model as needed. If the quantities $\hat{E}_i, \; i=1,\ldots,m$ are e.g. $l_p$ norms of the Fourier modes, then we can multiply Equations \eqref{reduced} with appropriate quantities and combine with Equations \eqref{conditions} to get
\begin{align*}\label{reduced}
\frac{d\hat{E}_1(\hat{u})}{dt} &= \sum_{i=1}^{m} a^{(1)}_i B_{1i} (t,\hat{u}(t)) \\
\frac{d\hat{E}_2(\hat{u})}{dt} &= \sum_{i=1}^{m} a^{(1)}_i B_{2i} (t,\hat{u}(t)) \\
 \quad     \cdots       \quad       & =   \quad       \cdots \quad \\
\frac{d\hat{E}_m(\hat{u})}{dt} &= \sum_{i=1}^{m} a^{(1)}_i B_{mi} (t,\hat{u}(t)) 
\end{align*}
where $B_{ij}=\p{}{a^{(1)}_j}\biggl(\frac{d\hat{E}_i(\hat{u}')}{dt}\biggr), \; i,j=1,\ldots,m$ are the new RHS functions that appear. Note that the RHS of the equations above does not involve primed quantities. The reason is that here the reduced quantities are computed by using the values of the resolved modes from the full system. The above system of equations is a linear system of equations for the vector of coefficients $a^{(1)}.$ The linear system can be written as
\begin{equation}\label{alphasystem}
B a^{(1)}= {\bf e}
\end{equation}
where ${\bf e}=\bigl(\frac{d\hat{E}_1(\hat{u})}{dt}, \ldots, \frac{d\hat{E}_m(\hat{u})}{dt} \bigr).$ This system of equations can provide us with the time evolution of the vector $a^{(1)}.$ 

The determination of coefficients for the reduced model through the system \eqref{alphasystem} is a time-dependent version of the method of moments. We specify the coefficients of the reduced model so that the reduced model reproduces the rates of change of a finite number of moments of the solution of the original system. This ensures that each term in the model is properly weighted so that the resulting reduced model reproduces, at the scales accessible to the reduced model, the dynamics (see \eqref{conditions}) of the original system. 

By construction, the entry $B_{ij}, \; i,j=1,\ldots,m$ of the matrix $B$ measures the contribution of the $j$-th term of the reduced model to the rate of change of $\hat{E}_i.$ In fact, the $j$-th column of the matrix $B$ is comprised of all the contributions of the $j$-th term in the reduced model to the rates of change of the different $\hat{E}_i.$ While the reduced system has no need to transfer activity from the resolved to the unresolved scales, the columns of $B$ corresponding to the activity-transferring terms will be zero (to the numerical precision used). Thus, the matrix $B$ will be singular. This can be monitored by estimating the rank of the matrix through the Singular Value Decomposition (SVD) \cite{golub}. When the smallest singular value becomes non-zero for the numerical precision used the reduced system starts transferring activity to the unresolved scales. After that instant we can use the system \eqref{alphasystem} to estimate the coefficient vector $a^{(1)}$ (see \cite {S09} for more details). 

The accurate determination of the coefficients guarantees that the rate at which the reduced model transfers activity to the unresolved scales is in agreement with the rate of transfer dictated by the full system. As we will see in Section \ref{numerics}, the accurate determination of the coefficients is crucial and moreover it can be hard. This is because the matrix $B$ can turn out to be extremely ill-conditioned. Such ill-conditioning is to be expected because the contributions of the different terms in the reduced model to the rate of change of the quantities $\hat{E}_i$ can vary dramatically from term to term. In fact, it is this significant variation of the contributions of the various terms which gauges the importance of the different terms in the reduced model. Eventually, the difference of the importance of the different terms is reflected on the values of the coefficients of the terms. This will become apparent when we present the numerical results in Section \ref{numerics}.

\subsection{The Mori-Zwanzig formalism}\label{mzformalism} 

We have presented in the previous section an abstract way of writing the reduced system which does not make any reference to a specific method for obtaining the functions $R_k^{(1)} (t,\hat{u}'(t))\; k\in F $ appearing on the RHS of \eqref{reduced}.  In order to proceed we need to specify the functions $R_k^{(1)} (t,\hat{u}'(t)).$ We will do that through the Mori-Zwanzig formalism \cite{CHK00,CHK3}. 

Suppose we are given the full system 
\begin{equation}\label{odes}
\frac{du(t)}{dt} = R (t,u(t)),
\end{equation}
where $u = ( \{u_k\}), \; k \in F \cup G$ 
with initial condition $u(0)=u_0.$ The system of ordinary differential equations
we are asked to solve can be transformed into a system of  linear
partial differential equations
\begin{equation}
\label{pde}
\pd{\phi_k}{t}=L \phi_k, \qquad \phi_k (u_0,0)=u_{0k}, \, k \in F \cup G
\end{equation}
where $L=\sum_{k \in F \cup G } R_i(u_0) \frac{\partial}{\partial u_{0i}}.$ The solution of \eqref{pde} is
given by $u_k (u_0,t)=\phi_k(u_0,t)$. Using semigroup notation we can rewrite (\ref{pde}) as
$$\pd{}{t} e^{tL} u_{0k}=L e^{tL} u_{0k}$$
Suppose that the vector of initial conditions can be divided as $u_0=(\hat{u}_0,\tilde{u}_0),$ where 
$\hat{u}_0$ is the vector of the resolved variables and $\tilde{u}_0$ is the vector of the unresolved variables.  Let $P$ be an orthogonal projection on the space of functions of $\hat{u}_0$ and $Q=I-P.$ We delay the specification of the projection operator until Section \ref{numerics} where we present applications to specific PDEs.

Equation \eqref{pde} 
can be rewritten as 
\begin{equation}
\label{mz}
\frac{\partial}{\partial{t}} e^{tL}u_{0k}=
e^{tL}PLu_{0k}+e^{tQL}QLu_{0k}+
\int_0^t e^{(t-s)L}PLe^{sQL}QLu_{0k}ds, \, k \in F,
\end{equation}
where we have used Dyson's formula
\begin{equation}
\label{dyson1}
e^{tL}=e^{tQL}+\int_0^t e^{(t-s)L}PLe^{sQL}ds.
\end{equation}
Equation (\ref{mz}) is the Mori-Zwanzig identity. 
Note that
this relation is exact and is an alternative way
of writing the original PDE. It is the starting
point of our approximations. Of course, we
have one such equation for each of the resolved
variables $u_k, k \in F$. The first term in (\ref{mz}) is
usually called Markovian since it depends only on the values of the variables
at the current instant, the second is called "noise" and the third "memory". 

If we write
$$e^{tQL}QLu_{0k}=w_k,$$ 
$w_k(u_0,t)$ satisfies the equation
\begin{equation}
\label{ortho}
\begin{cases}
&\frac{\partial}{\partial{t}}w_k(u_0,t)=QLw_k(u_0,t) \\ 
& w_k(u_0,0) = QLx_k=R_k(u_0)-(PR_k)(\hat{u_0}). 
\end{cases} 
\end{equation}
If we project (\ref{ortho}) we get
$$P\frac{\partial}{\partial{t}}w_k(u_0,t)=
PQLw_k(u_0,t)=0,$$
since $PQ=0$. Also for the initial condition
$$Pw_k(u_0,0)=PQLu_{0k}=0$$
by the same argument. Thus, the solution
of (\ref{ortho}) is at all times orthogonal
to the range of $P.$ We call
(\ref{ortho}) the orthogonal dynamics equation. Since the solutions of 
the orthogonal dynamics equation remain orthogonal to the range of $P$, 
we can project the Mori-Zwanzig equation (\ref{mz}) and find
\begin{equation}
\label{mzp}
\frac{\partial}{\partial{t}} Pe^{tL}u_{0k}=
Pe^{tL}PLu_{0k}+
P\int_0^t e^{(t-s)L}PLe^{sQL}QLu_{0k} ds.
\end{equation}
In order to proceed with the computation of the reduced model we need to compute the Markovian term and the memory term. While the Markovian term is usually rather straightforward to compute (see Section \ref{numerics}), the memory term computation is rather involved due to the presence of the orthogonal dynamics evolution operator $e^{tQL}.$ In fact, it is the presence of this operator which makes, in general, the computation of MZ reduced models prohibitively expensive (see \cite{CS05} for a thorough discussion). 

One can start from \eqref{mzp} and based on assumptions derive simplified reduced models that are easier to calculate \cite{CHK3,bernstein,HS06,S06}. However, it has proven difficult to justify the assumptions underlying the derivation of the simplified models. In addition, it is hard to work with the expression for the memory term given in \eqref{mzp}. We offer here an alternative presentation of the memory term which allows us to fit the MZ formalism in the framework described in Section \ref{full_reduced}.

Note that the memory term $\int_0^t e^{(t-s)L}PLe^{sQL}QLu_{0k}$ can be written, through Dyson's formula \eqref{dyson1}, as 
\begin{equation}\label{dyson2}
\int_0^t e^{(t-s)L}PLe^{sQL}QLu_{0k} ds =e^{tL}QLu_{0k} -e^{tQL}QLu_{0k} 
\end{equation}
which can be rewritten (using the linearity of the operator $e^{tL}$) as
\begin{equation}\label{dyson3}
e^{tL}QLu_{0k} -e^{tQL}QLu_{0k} = e^{tL} \bigl( QLu_{0k} -e^{-tL} e^{tQL}QLu_{0k} \bigr).
\end{equation}
Since $I=P+Q$ we find
\begin{equation}\label{dyson4}
\int_0^t e^{(t-s)L}PLe^{sQL}QLu_{0k} ds=e^{tL} \bigl( QLu_{0k} -e^{C(t,u_0)}QLu_{0k} \bigr)
\end{equation}
where $C(t,u_0)=-tPL+ [-tL,tQL] + \ldots $ is given by the BCH series and $[-tL,tQL]=-tLtQL-tQL(-tL)$ (see e.g \cite{bellman}). We have $[-tL,tQL]=[tL,tPL]=[tQL,tPL].$ In general, $[tQL,tPL] \neq 0.$ However, depending on the initial conditions, $[tQL,tPL]$ may be small and thus allow the simplification of the memory term expression. The magnitude of $[tQL,tPL]$ will be examined in detail in a future publication. 

If we assume that $C(t,u_0)\approx -tPL,$ expansion of the operator $e^{-tPL}$ in Taylor series around $t=0$ gives 
\begin{gather}\label{dyson6}
P\int_0^t e^{(t-s)L}PLe^{sQL}QLu_{0k} ds \approx \\
 \sum_{j=1}^{\infty} (-1)^{j+1} \frac{t^j}{j!} Pe^{tL}(PL)^jQLu_{0k}. \notag
\end{gather}
One can obtain different simplified models by truncating the series in \eqref{dyson6} after different values of $j.$ In particular, if we omit all the terms after the first one we obtain the $t$-model which has been studied thoroughly \cite{CHK3,bernstein,HS06}.

%Since $I=P+Q,$ the second term in the parenthesis can be simplified and we have   
%\begin{equation}\label{dyson4}
%e^{tL}QLu_{0k} -e^{tQL}QLu_{0k} = e^{tL} \bigl( QLu_{0k} -e^{-tPL}QLu_{0k} \bigr).
%\end{equation}
%Now, we expand the operator $e^{-tPL}$ in Taylor series around $t=0$ and we find for the memory term
%\begin{equation}\label{dyson5}
%\int_0^t e^{(t-s)L}PLe^{sQL}QLu_{0k} ds =  \sum_{j=1}^{\infty} (-1)^{j+1} \frac{t^j}{j!} e^{tL}(PL)^jQLu_{0k}.
%\end{equation}
%Finally, for the projected memory term $P\int_0^t e^{(t-s)L}PLe^{sQL}QLu_{0k} ds$ we find
%\begin{equation}\label{dyson6}
%P\int_0^t e^{(t-s)L}PLe^{sQL}QLu_{0k} ds =  \sum_{j=1}^{\infty} (-1)^{j+1} \frac{t^j}{j!} Pe^{tL}(PL)^jQLu_{0k}.
%\end{equation}
%One can obtain different simplified models by truncating the series in \eqref{dyson6} after different values of $j.$ In particular, if we omit all the terms after the first one we obtain the $t$-model which has been studied thoroughly \cite{CHK3,bernstein,HS06}.

The series representation of the memory term in \eqref{dyson6} is based on the assumption of analyticity in time of the operator $e^{-tPL}.$ This assumption may be true for small $t$ but it does not have to hold for larger $t.$ In other words, the Taylor expansion of the operator $e^{-tPL}$  has, in general, only a {\it finite} radius of convergence. Insisting on using the Taylor expansion of the operator $e^{-tPL}$  as is for later times is dangerous and can lead to the instability of the reduced model (see also Section \ref{renovsnoreno}). In fact, when dealing with full systems coming from discretizations of singular PDEs, the breakdown of the Taylor expansion of the operator $e^{-tPL}$  is related to the onset of underresolution on the part of the full system. 

To proceed we need to put the MZ model given by \eqref{mzp} and \eqref{dyson6} in the framework of Section \ref{full_reduced}. To do that we set
\begin{align}
 R^{(1)}_{1k} &= Pe^{tL}PLu_{0k}, \\
 R^{(1)}_{jk} &=  (-1)^{j} \frac{t^{j-1}}{(j-1)!} Pe^{tL}(PL)^{j-1}QLu_{0k},  \;\; j=2,\ldots.
\end{align}
With this identification we have, in essence, embedded the reduced models derived through the MZ formalism in a larger class of reduced models which share the same functional form with the MZ models but who are allowed to have different coefficients. In the notation of Section \ref{full_reduced}, the original MZ models correspond to the coefficient vector $a^{(1)}=(1,1,1,...).$ 

While the original MZ models may suffer from instabilities, the new models can be made stable by {\it assigning to each term in the reduced model the appropriate coefficient}. The magnitude of the coefficient of a term reflects the importance of the term in the reduced model. The values of the coefficients can now be determined by solving the linear algebraic system \eqref{alphasystem}. This ensures that the coefficient of each term in the model is properly redefined (renormalized) so that the resulting reduced model reproduces, at the scales accessible to the reduced model, the dynamics (see \eqref{conditions}) of the original system. 

\subsection{The renormalized Mori-Zwanzig (rMZ) algorithm}\label{combine}

We are now in a position to state the renormalized Mori-Zwanzig algorithm which constructs a reduced model with the necessary coefficients computed on the fly. 

\vskip14pt
{\bf Renormalized Mori-Zwanzig Algorithm}
\begin{enumerate}
\item
Choose a number of terms, say $m$, to keep at the Taylor expansion of the memory term. 
\item
Evolve the full system and compute, at every step, using the SVD, the rank of the $(m+1)\times (m+1)$ matrix $B.$
\item
When the smallest singular value $\sigma_{m+1}=TOL,$ (we assume that the singular values are indexed from largest to smallest) where $TOL$ a prescribed tolerance, solve the system \eqref{alphasystem} for the coefficients.
\item
For the remaining simulation time, switch from the full system to the reduced model with the estimated values of the coefficients.  
\end{enumerate}
In order to apply the algorithm above, we need to specify the quantities $\hat{E}_i, \; i=1,\ldots,m.$ Also, we need to compute the expression for the Markovian term, as well as the expressions for the terms in the Taylor expansion of the memory term. This can be a formidable task. However, as we will show in the next section, for the case of Burgers and Euler, these computations can be automated and thus we can compute arbitrarily high order terms in the memory term expansion with minimal effort and storage.

%%%%%%%%%%%%%%End of Section Renormalized MZ%%%%%%%%%%%%

\section{Application of rMZ to 1D Burgers and 3D Euler equations}\label{numerics}

In this section we present results of the rMZ algorithm for the 1D Burgers and the 3D Euler equations.

\subsection{1D Burgers equation}\label{burgers}

\subsubsection{Setup of the reduced model}\label{burgers-setup}
We use the 1D inviscid Burgers equation as an instructive example for the constructions presented in this section. The equation is given by 
\begin{equation}\label{burgersequation}
u_t+u u_x = 0.
\end{equation}
Equation (\ref{burgersequation}) should be supplemented with an initial condition $u(x,0)=u_0(x)$ and boundary conditions. We solve (\ref{burgersequation}) in the interval $[0,2\pi]$ with periodic boundary conditions. This allows us to expand the solution in Fourier series
$$u_{M}(x,t )=\underset{k \in F \cup G}{\sum} u_k(t) e^{ikx},$$
where $F \cup G=[-\frac{M}{2},\frac{M}{2}-1].$ We have written the set of Fourier modes as the union of two sets 
in anticipation of the construction of the reduced model comprising only of the modes in $F=[-\frac{N}{2},\frac{N}{2}-1],$ where $ N < M.$
The equation of motion for the Fourier mode $u_k$ becomes
\begin{equation}
\label{burgersode}
 \frac{d u_k}{dt}=- \frac{ik}{2} \underset{p, q \in F \cup G}{\underset{p+q=k  }{ \sum}} u_{p} u_{q}.
\end{equation}
To conform with 
the Mori-Zwanzig formalism we set 
$$R_k(u)=- \frac{ik}{2} \underset{p, q \in F \cup G}{\underset{p+q=k  }{ \sum}}u_{p} u_{q}$$
and we have 
\begin{equation}
\label{burgersodemz}
\frac{d u_k}{dt}=R_k(u) 
\end{equation}
for $ k \in F \cup G.$ 
The system (\ref{burgersodemz}) is supplemented by the initial 
condition $u_0=(\hat{u}_0,\tilde{u}_0)=(\hat{u}_0,0).$ We focus on initial conditions where 
the unresolved Fourier modes are set to zero. We also define $L$ by 
$$L=\sum_{k \in F \cup G} R_k(u_0) \frac{\partial}{\partial u_{0k}}.$$ 
Note that $Lu_{0k}=R_k(u_0).$ 

We also need to define a projection operator $P.$ For a function $h(u_0)$ of all the 
variables, the projection operator we will use is defined by $P(h(u))=P(h(\hat{u}_0,\tilde{u}_0))=h(\hat{u}_0,0),$ i.e. 
it replaces the value of the unresolved variables $\tilde{u}_0$ in any function $h(u_0)$ by zero. Note that this choice of projection is consistent with the initial conditions we have chosen. Also, we define the Markovian term 
$$\hat{R}^{(1)}_1k(\hat{u}_0)= PLu_{0k}=PR_k(u_0)=- \frac{ik}{2} \underset{p, q \in F}{\underset{p+q=k  }{ \sum}}  \hat{u}_{0p}  \hat{u}_{0q}.$$ 
The Markovian term has the same functional form as the RHS of the full system but is restricted to a sum over only the resolved modes in $F.$ The full system conserves the energy $\frac{1}{2}\sum_{k \in F \cup G} |u_k|^2$ contained in all the modes. Similarly, the Markovian term of the reduced model does {\it not} alter the energy content of the resolved modes. The necessary energy transfer out of the resolved modes rests on the memory terms. Based on our choice of projection operator and the scaling symmetries of the Burgers equation we set $N=\frac{M}{2}.$

With the definition of $P$ given above, we find for $QLu_{0k}$
$$QLu_{0k}= - \frac{ik}{2} \underset{q \in G}{\underset{p \in F}{\underset{p+q=k  }{ \sum}}}  u_{0p}  u_{0q}  - \frac{ik}{2} \underset{q \in F}{\underset{p \in G}{\underset{p+q=k  }{ \sum}}}  u_{0p}  u_{0q}  
- \frac{ik}{2} \underset{q \in G}{\underset{p \in G}{\underset{p+q=k  }{ \sum}}}  u_{0p}  u_{0q}.$$
The expression for $QLu_{0k}$ contains three terms which involve at least one wavenumber in the unresolved range $G.$ The terms in the Taylor expansion of the memory term are given by 
$$R^{(1)}_{jk} =  (-1)^{j} \frac{t^{j-1}}{(j-1)!} Pe^{tL}(PL)^{j-1}QLu_{0k},  \;\; j=2,\ldots $$ 
For the $j$-th term we have 
\begin{equation}\label{memoryterm}
(PL)^{j-1}QLu_{0k}=(PL)^{j-1} \biggl(   - \frac{ik}{2} \underset{q \in G}{\underset{p \in F}{\underset{p+q=k  }{ \sum}}}  u_{0p}  u_{0q}  - \frac{ik}{2} \underset{q \in F}{\underset{p \in G}{\underset{p+q=k  }{ \sum}}}  u_{0p}  u_{0q}  - \frac{ik}{2} \underset{q \in G}{\underset{p \in G}{\underset{p+q=k  }{ \sum}}}  u_{0p}  u_{0q} \biggr) .
\end{equation}

\subsubsection{Recursive computation of the memory terms}\label{burgers-recursive}
The expression for $(PL)^{j-1}QLu_{0k}$ in \eqref{memoryterm} for the $j$th term ($j=2,\ldots$) can be computed recursively using a simple construction based on a Pascal triangle. Note that for our choice of projection operator $P,$ we have 
$$(PL)^{j-1}  \biggl( - \frac{ik}{2} \underset{q \in G}{\underset{p \in G}{\underset{p+q=k  }{ \sum}}}  u_{0p}  u_{0q} \biggr) =0 .$$
We begin with the (first-order) term for $j=2$ i.e., $PLQLu_{0k}.$ We find
\begin{equation}\label{firstorder}
PLQLu_{0k}=   -2 \frac{ik}{2} \underset{q \in G}{\underset{p \in F}{\underset{p+q=k  }{ \sum}}}  Pu_{0p}  PLu_{0q}. 
\end{equation}
The first order term can be computed by convolving the resolved part of $Pu_{0p}$ with the unresolved part of $PLu_{0q}.$ In practice, all the convolutions sums can be computed using Fast Fourier Transforms \cite{boyd}. Note that the expression $Pu_{0p}$ is linear in the Fourier modes while $PLu_{0q}$ is quadratic. Thus, the convolution sum in $PLQLu_{0k}$ (including the factor $-\frac{ik}{2}$)can be denoted as $(1r*2u),$ where $*$ stands for convolution while $r$ and $u$ stand for the resolved and unresolved parts. This notation facilitates the recognition of the pattern for the higher order terms. With this notation, the first-order term can be written as 
\begin{equation}\label{firstorder2}
PLQLu_{0k}=  2 \times {\bf 1} (1r * 2u)
\end{equation} 
where we have used bold face to denote the coefficient. We continue with the second order term $PLPLQLu_{0k}.$ We find
\begin{equation}\label{secondorder}
PLPLQLu_{0k}=  2 \biggl( - \frac{ik}{2} \underset{q \in G}{\underset{p \in F}{\underset{p+q=k  }{ \sum}}}  Pu_{0p}  PLPLu_{0q} -  \frac{ik}{2} \underset{q \in G}{\underset{p \in F}{\underset{p+q=k  }{ \sum}}}  PLPu_{0p}  PLu_{0q} \biggr).  
\end{equation}
The convolution sums in this term can be denoted as $(1r * 3u)$ and $(2r * 2u).$ The second order term can be written as
\begin{equation}\label{secondorder2}
PLPLQLu_{0k}=  2 \times  \biggl( {\bf 1} (1r * 3u) + {\bf 1} (2r * 2u)  \biggr)
\end{equation} 
To see the pattern more clearly we need one more term,  $PLPLPLQLu_{0k}.$ We find
\begin{gather}\label{thirdorder}
PLPLPLQLu_{0k}= \\
2 \biggl(  - \frac{ik}{2} \underset{q \in G}{\underset{p \in F}{\underset{p+q=k  }{ \sum}}}  Pu_{0p}  PLPLPLu_{0q}  \notag \\
 - 2 \frac{ik}{2} \underset{q \in G}{\underset{p \in F}{\underset{p+q=k  }{ \sum}}}  PLPu_{0p}  PLPLu_{0q} \notag  \\
 - \frac{ik}{2} \underset{q \in G}{\underset{p \in F}{\underset{p+q=k  }{ \sum}}}  PLPLPu_{0p}  PLu_{0q} \biggr) \notag
\end{gather}
The terms in the parenthesis can be denoted as $(1r * 4u),$ $(2r * 3u)$ and $(3r * 2u).$ The third order term can be written as
\begin{equation}\label{thirdorder2}
PLPLPLQLu_{0k}=  2 \times \biggl( {\bf 1} (1r * 4u) + {\bf 2} (2r * 3u) + {\bf 1} (3r * 2u)  \biggr)
\end{equation} 
By examining the expressions in \eqref{firstorder2}-\eqref{thirdorder2} we see that the memory terms can be computed as weighted sums of convolution sums where the weights are given by appropriate Pascal triangle coefficients (the bold face numbers). This was to be expected since we started with a convolution sum (involving products) of two functions and each new term in the Taylor series involves a differentiation . Moreover, the number of convolution sums that need to be added is equal to the order of the memory term in the Taylor expansion. Also, for each term, the convolution sums involve expressions whose degree (in Fourier modes) follows an easily discernible pattern. For the $l$th order term in the Taylor series we need the convolution sums $(1r*(l+1)u), (2r*lu), \ldots, (lr * 2u).$ 

Finally, the expressions entering the convolution sums can also be computed by a Pascal triangle construction.  For example, in order to calculate the third order term, as can be seen from \eqref{thirdorder}, one needs to compute and store only the quantities $Pu_{0p}, PLu_{0p}, PLPLu_{0p}, PLPLPLu_{0p}$ for $p \in F \cup G.$ Note that the quantities $Pu_{0p}, PLPu_{0p}$ and $PLPLPu_{0p}$ which are also needed are the same as  $Pu_{0p}, PLu_{0p}, PLPLu_{0p}$ for the {\it resolved} modes and zero for the {\it unresolved} modes. So, they do not need to be stored. They can be quickly constructed when needed. We see that the storage requirements for the calculation of the memory terms grows only linearly in the order of the Taylor expansion. Also, the ability to calculate the needed expressions through FFTs speeds up significantly the calculation of the various memory terms.

The recursive estimation of the memory terms allows us to calculate memory terms of very high order efficiently, without having to write down explicitly the analytical expressions which become very complicated after the first few orders in the expansion.

\subsubsection{Results using rMZ}\label{burgers-numerics}

\begin{figure}
\centering
\epsfig{file=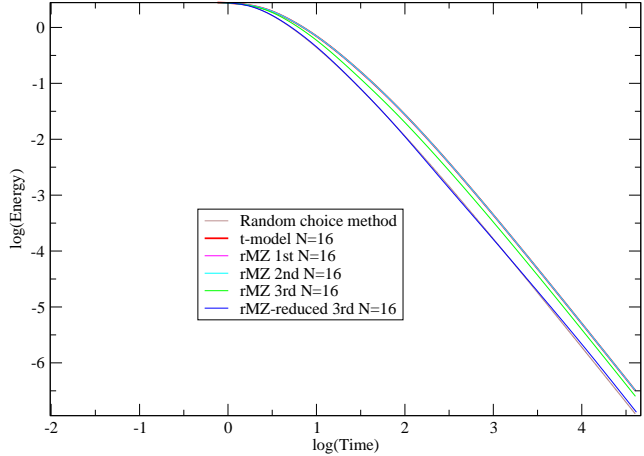,width=4.in}
\caption{1D Burgers equation. Comparison of the evolution of the energy in the resolved modes computed by the random choice method, the $t$-model and various rMZ models (see text for details).}
\label{plot_burgers_collective}
\end{figure}

Figure \ref{plot_burgers_collective} shows the evolution of the energy $\frac{1}{2}\sum_{k \in F} |u_k|^2$ of the resolved modes computed by reduced models of different orders and the random choice method \cite{chorinran}. The initial condition is $u_0(x)=\sin x.$ This leads to the formation of a standing shock at $T=1.$ All the reduced models use $N=16$ Fourier modes while the full system has $M=32$ modes. The results of the reduced models are compared to a converged solution of the random choice method with $N=4096$ points. The energy of the random choice method solution was computed using only $N=16$ modes. However, note that practically all the energy of the random choice method solution is concentrated in the first few Fourier modes, so even if we had computed the energy for all $N=4096$ Fourier methods the results would not have changed. This is to be expected, since for the initial condition we are using, a standing shock forms at time $T=1$ and, thus, by time $T=100$ the only Fourier modes having some energy left in them are the first few.   

The quantities $\hat{E}_i$ used to set up the linear algebraic system needed to compute the coefficients of the reduced system are $l_p$ norms of the solution. In particular, for the first order model we use $\hat{E}_i=\sum_{k \in F} |u_k|^{2i}, \; i=1,2.$  For the second order model $\hat{E}_i=\sum_{k \in F} |u_k|^{2i}, \; i=1,2,3$ and for the third order model $\hat{E}_i=\sum_{k \in F} |u_k|^{2i}, \; i=1,2,3,4.$ In general, for the reduced model of order $\lambda$ we need $\lambda+1$ quantities because we also have to compute the coefficient of the Markovian term. All the calculations are done in double precision. The tolerance $TOL$ which is used to decide when it is time to switch to the reduced model is set to $TOL=10^{-12}.$ The systems of ordinary differential equations for the different reduced models were solved using the Runge-Kutta-Fehlberg method with the stepsize control tolerance set to $10^{-10}$ \cite{hairer}.

The numerical problem of solving the linear system for the coefficients is hard because the resulting system has very large condition number and very small determinant. This happens for three reasons. First, the dominant contribution to the linear system matrix comes from the Markovian term (except for the contribution to the rate of change of the $\hat{E}_1=\sum_{k \in F} |u_k|^{2}$ which is zero). This means that the coefficient of the Markovian term is practically 1. Second, the contributions of each memory term to the rates of change of the different $\hat{E}_i$ vary dramatically. Third, the contributions to each $\hat{E}_i$ by the different memory terms also varies substantially. Of course, this situation is exacerbated if we use more terms we use in the expansion. For the case when we retain up to the third order term in the memory expansion, we have to deal with condition numbers of the order $10^{11}$ and determinant values of order $10^{-20}.$ Inevitably, even the use of double precision cannot provide us with an accurate estimate of the coefficients. Since the linear system matrix is practically singular (for the numerical precision used) we have chosen to solve the linear system using the SVD algorithm \cite{golub}. 

A partial remedy to the problem comes from a slight modification in the way of estimating the coefficients. Since we know that the Markovian term coefficient is practically 1, we can set it to 1, and subtract the column of contributions of the Markovian term from the RHS of the linear system. This allows us to reduce the dimensionality of the linear system to be solved from $(\lambda+1) \times (\lambda+1)$ to $\lambda \times \lambda.$ This practice of subtracting almost equal numbers is not advisable in general because it leads to loss of significant digits \cite{golub}. However, in our case it helps to improve the results by lowering the condition number of the matrix from about $10^{11}$ to about $10^{5}.$ In Figure 1, the estimation of the coefficients for the third order model using the reduced dimension matrix is denoted as "rMZ-reduced 3rd".

As shown in Figure \ref{plot_burgers_collective}, the rMZ models of first and second order give practically the same results as the $t$-model. The third order model gives a slight improvement. However, when the reduced dimension matrix is used, the energy evolution predicted by the third order model is practically identical to the correct energy evolution of the resolved modes predicted by the random choice method. If we increase the resolution of the reduced model, the numerical problems for the calculation of higher order coefficients become even more pronounced. This is to be expected, since a larger resolution means that the renormalized coefficients of the reduced model will be smaller. Thus, computing them with accuracy is more difficult. 

We have to note, that the computation of the higher order coefficients is more difficult for our choice of initial condition since it only involves one active Fourier mode. Initial conditions with more active Fourier modes will transfer activity to the unresolved scales at a higher rate and thus the corresponding renormalized coefficients will be larger. A detailed study of the behavior of the coefficients for different initial conditions will be presented elsewhere. For the first order model, we have already presented in \cite{S09} a detailed study about the change of the value of the renormalized coefficient with resolution up to the order of $10^5$ Fourier modes. In that work, the renormalized coefficient calculation was used to determine, in a fixed point analysis, the blow-up exponent.

%%%%%%%%%%%%%%End of Section Burgers%%%%%%%%%%%%%%%%%

\subsection{3D incompressible Euler equations}\label{euler}

Consider the 
3D incompressible Euler equations with periodic boundary conditions in the cube $[0,2\pi]^3$:

\begin{gather}
\label{eulerpde}
u_t+u\cdot \nabla u= - \nabla p  ,\; \nabla \cdot u=0, 
\end{gather}
where $u(x,t)=(u_1(x_1,x_2,x_3,t),u_2(x_1,x_2,x_3,t),u_3(x_1,x_2,x_3,t))$ is the velocity, $p$ is the 
pressure and $\nabla= (\frac{\partial}{\partial x_1},\frac{\partial}{\partial x_2},\frac{\partial}{\partial x_3} ).$ The system in (\ref{euler}) is 
supplemented with the initial condition $u(x,0)=u_0(x)$ which is also periodic and incompressible and 
$x=(x_1,x_2,x_3).$ Since we are working with periodic boundary conditions, we expand the solution in Fourier series 
keeping $N$ modes in each spatial direction, 
$$u_{M}(x,t )=\underset{k \in F \cup G}{\sum} u_k(t) e^{ikx},$$
where $F \cup G=[-\frac{M}{2},\frac{M}{2}-1]\times[-\frac{M}{2},\frac{M}{2}-1]\times[-\frac{M}{2},\frac{M}{2}-1].$ Also $k=(k_1,k_2,k_3)$ and $u_k(t)=(u^1_k(t),u^2_k(t),u^3_k(t)).$ 

The equation of motion for the Fourier mode $u_k$ becomes
\begin{equation}
\label{eulerode}
 \frac{d u_k}{dt}=- i \underset{p, q \in F \cup G}{\underset{p+q=k  }{ \sum}} k \cdot u_{p} 
A_{k} u_{q},
\end{equation}
where $A_k= I - \frac{k k^T}{|k|^2}$ is the incompressibility projection matrix and $I$ is the $3\times3$ 
identity matrix. The symbol $\cdot$ denotes inner product in $\mathbb{R}^3.$ The system (\ref{eulerode}) is supplemented by the initial condition $u_0=\{u_k(0)\}=\{u_{0k}\}, \; k \in F \cup G,$ where $u_{0k}$ are the Fourier coefficients of the initial condition $u_0(x).$ 

To conform with 
the MZ formalism we set 
$$R_k(u)=- i \underset{p, q \in F \cup G}{\underset{p+q=k  }{ \sum}} k \cdot u_{p} A_{k} u_{q}$$
and we have 
\begin{equation}
\label{eulerodemz}
\frac{d u_k}{dt}=R_k(u) 
\end{equation}
for $ k \in F \cup G.$ 
The system (\ref{eulerodemz}) is supplemented by the initial 
condition $u_0=(\hat{u}_0,\tilde{u}_0)=(\hat{u}_0,0).$ Note that we focus on initial conditions where 
the unresolved Fourier modes are set to zero. We also define $L^0$ by 
$$L=\sum_{k \in F \cup G} R_k(u_0) \frac{\partial}{\partial u_{0k}}.$$ 
Note that $Lu_{0k}=R_k(u_0).$ Consider the subset $F =[-\frac{N}{2},\frac{N}{2}-1]\times[-\frac{N}{2},\frac{N}{2}-1]\times[-\frac{N}{2},\frac{N}{2}-1]$ for $N < M.$ We will construct the reduced models for the Fourier modes $u_k$ with $k \in F. $ 

We need to define a projection operator $P.$ For a function $h(u_0)$ of all the 
variables, the projection operator we will use is defined by $P(h(u))=P(h(\hat{u}_0,\tilde{u}_0))=h(\hat{u}_0,0),$ i.e. 
it replaces the value of the unresolved variables $\tilde{u}_0$ in any function $h(u_0)$ by zero. Note, that this choice of projection is consistent with the initial conditions we have chosen. Based on our choice of projection operator and the scaling symmetries of the Euler equations we set $N=\frac{M}{2}.$ 

Define  
$$\hat{R}_k(\hat{u}_0)= PR_k(u_0)=- i \underset{p, q \in F}{\underset{p+q=k  }{ \sum}} k \cdot \hat{u}_{0p} A_{k} \hat{u}_{0q}.$$
The Markovian term has the same functional form as the RHS of the full system but is restricted to a sum over only the resolved modes in $F.$ The full system conserves the energy $\frac{1}{2}\sum_{k \in F \cup G} |u_k|^2$ contained in all the modes. Similarly, the Markovian term of the reduced model does {\it not} alter the energy content of the resolved modes. The necessary energy transfer out of the resolved modes rests on the memory terms. 

With the definition of $P$ given above, we find for $QLu_{0k}$
$$QLu_{0k}= - i \underset{q \in G}{\underset{p \in F}{\underset{p+q=k  }{ \sum}}} k \cdot u_{0p} 
A_{k} u_{0q}  -  \underset{q \in F}{\underset{p \in G}{\underset{p+q=k  }{ \sum}}} k \cdot u_{0p} 
A_{k} u_{0q}  
- i \underset{q \in G}{\underset{p \in G}{\underset{p+q=k  }{ \sum}}}  k \cdot u_{0p} 
A_{k} u_{0q}.$$
The expression for $QLu_{0k}$ contains three terms which involve at least one wavenumber in the unresolved range $G.$ The terms in the Taylor expansion of the memory term are given by 
$$R^{(1)}_{jk} =  (-1)^{j} \frac{t^{j-1}}{(j-1)!} Pe^{tL}(PL)^{j-1}QLu_{0k},  \;\; j=2,\ldots $$ 
For the $j$-th term we have 
\begin{gather}\label{memoryterm2}
(PL)^{j-1}QLu_{0k}=(PL)^{j-1} \biggl(   - i \underset{q \in G}{\underset{p \in F}{\underset{p+q=k  }{ \sum}}}  k \cdot u_{0p} 
A_{k} u_{0q}  - i \underset{q \in F}{\underset{p \in G}{\underset{p+q=k  }{ \sum}}}  k \cdot u_{0p} 
A_{k} u_{0q}  \\
- i \underset{q \in G}{\underset{p \in G}{\underset{p+q=k  }{ \sum}}}  k \cdot u_{0p} 
A_{k} u_{0q} \biggr) .
\end{gather}
The different terms in the memory expansion can be computed recursively as in the case of Burgers. However, there is a slight complication because the presence of the incompressibility operator and of the inner product on the RHS destroys the commutativity which allowed us in Burgers to group terms (the factor 2 which appears outside every parenthesis there). This problem can be addressed by a construction which uses 2 Pascal triangles instead of 1 used in the case of Burgers. For each order, one adds up the corresponding terms from the 2 Pascal triangles and obtains the desired memory term. Other than that, the recursive algorithm remains the same and we omit the details. Also, note that all the higher order terms are divergence-free by construction.

\begin{figure}
\centering
\epsfig{file=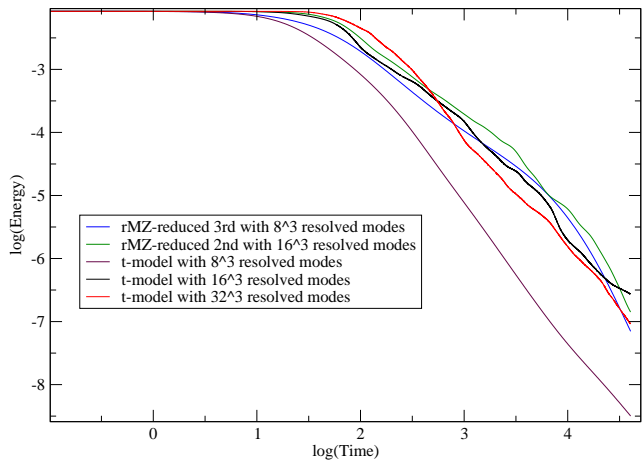,width=4.in}
\caption{3D Euler equation. Comparison of the evolution of the energy in the resolved modes computed by the $t$-model and various rMZ models (see text for details).}
\label{plot_euler_collective}
\end{figure}

We have used the same quantities $\hat{E}_i$ as in the case of Burgers, with the obvious generalizations, since for 3D Euler we have a 3-dimensional velocity vector instead of the scalar velocity in Burgers. Also, even though for 3D Euler we have a 3-dimensional vector, we have assumed that the reduced model renormalized coefficients are the same for all 3 velocities. This is a simplifying assumption. Of course, one can use different renormalized coefficients for the different velocities at the expense of having to solve a larger linear system for the renormalized coefficients. A detailed study of that case will be presented in future work.

We have used the Taylor-Green initial condition (see e.g. \cite{gottlieb} and references therein) which is given by 
\begin{eqnarray*}
u_1(x,0)&=& \sin(x_1)\cos(x_2)\cos(x_3), \\
u_2(x,0)&=&- \cos(x_1)\sin(x_2)\cos(x_3),\\
u_3(x,0)&=&  0.
\end{eqnarray*}

Figure \ref{plot_euler_collective} shows the evolution of the energy $\frac{1}{2}\sum_{k \in F} |u_k|^2$ of the resolved modes for different resolutions computed by rMZ reduced models of different orders and the $t$-model. We have presented results for the rMZ reduced models using the reduced linear system matrix approach discussed above to tame the condition number of the matrix. Based on these results, we make two observations. 

First, for $8^3$ resolved modes, the rMZ third order model dissipates energy at a slower rate than the $t$-model with $8^3$ resolved modes. This is true not only for the third order model but also for the first and second order models (we have omitted those results to avoid cluttering the figure). This slower rate of energy dissipation compared to the $t$-model holds also for the case of $16^3$ resolved modes. 

The second observation is that the rate of energy dissipation of the rMZ models is consistent with the rate predicted by the $t$-model with {\it higher} resolution. This is to be expected since a higher order model should result in a more accurate prediction of the energy dissipation rate. 

The reader may be concerned about the small resolutions used in the numerical experiments. There are two reasons for that. First, if one keeps several terms in the memory expansion, then, for a very smooth initial condition like the one we use, the matrix $B$ becomes even more ill-conditioned for large resolutions. However, this is not a severe problem. On the contrary, it signifies that most of the higher order terms should have small coefficients and thus can be safely removed from the model. 

The second reason we have used small resolutions both for Burgers and Euler is because an accurate reduced model should be able to reproduce the correct energy content for its resolved scales no matter how small the resolution. For example, for Burgers, where we know what the energy content should be after the singularity, we see that the rMZ model with a small resolution ($16^3$) indeed reproduces the correct energy content for this resolution. 

%%%%%%%%%%%%%End of Section Euler%%%%%%%%%%%%%%%%%

\subsection{Universality of the renormalized coefficients}\label{universality}

We have hinted in the introduction at the possibility that the renormalized coefficients may be determined by two factors: i) the ratio of the smallest active scale in the initial condition to the smallest resolvable scale and ii) the scaling symmetries of the equation under investigation. Even though the numerical difficulties with the ill-conditioned linear system matrix do not allow us at present to study accurately the higher order renormalized coefficients, we have enough accuracy to study the first order renormalized coefficient both for 1D Burgers and 3D Euler. Note that the two equations share the same scaling symmetries. Also, we have chosen for Burgers the initial condition $u_0(x)=\sin x$ which has only one active Fourier mode (for $k=\pm 1$) and for 3D Euler, the Taylor-Green initial condition which also has only active Fourier modes for $k_i = \pm 1, \; i=1,2,3.$ 

\begin{figure}
\centering
\epsfig{file=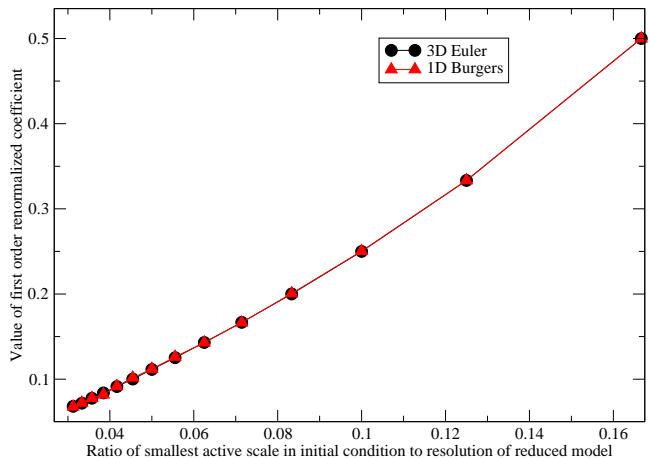,width=4.in}
\caption{Comparison, for 1D Burgers and 3D Euler,  of the value of the renormalized coefficient for the first order rMZ model for different resolutions. The initial conditions for Burgers and Euler have only 1 active Fourier mode in each direction (see text for details).}
\label{plot_burgers_euler_coefficient}
\end{figure}

Figure \ref{plot_burgers_euler_coefficient} shows the comparison of the value of the renormalized first-order coefficient for Burgers and 3D Euler as a function of the ratio of the smallest scale active in the the initial condition to the smallest scale of the reduced model.  We make two observations. First, the values of the renormalized coefficient for the two equations are in remarkable agreement. Second, from the slope of the linear fit, we see that the value of the coefficient is practically equal to the ratio of the smallest scale active in the the initial condition to the smallest scale of the reduced model. 

Needless to say that one example is not enough to infer the  generality of the result for arbitrary initial conditions. A theoretical explanation of this result is lacking at the moment. Note that due to the way we have defined the terms in the expansion of the memory, all the terms have the same dimensions as the Markovian term and the left-hand side of the equation for each Fourier mode. So, the corresponding coefficients have to be dimensionless. Thus, we expect the coefficients to depend on ratios of quantities with the same dimensions. Here we have investigated the possibility that this ratio is that of the smallest active scale in the initial condition to the smallest active scale of the reduced model.   

We should comment here on the behavior of the rMZ algorithm for the 2D Euler equations for which the 2D version of the Taylor-Green initial condition is an exact solution, i.e. a steady state. Exactly because it is a steady state there is no need for a reduced model. Application of the rMZ algorithm agrees with this. There is never any need to transfer energy to the unresolved scales and thus, no need to switch to a reduced model. The contributions of the different memory terms to the matrix $B$ are all well below the double precision threshold. This allows the freedom to assign to the renormalized coefficient the values shown in Figure \ref{plot_burgers_euler_coefficient} without incurring any trouble. In other words, the behavior of the solution of the 2D Euler for the Taylor-Green initial condition does not contradict the agreement for the renormalized coefficient of the 1D Burgers and 3D Euler equations shown in Figure \ref{plot_burgers_euler_coefficient}. 

%%%%%%%%%%%%%%End of Section Universality%%%%%%%%%%%%%%%%

\subsection{rMZ vs MZ}\label{renovsnoreno}

\begin{figure}
\centering
\epsfig{file=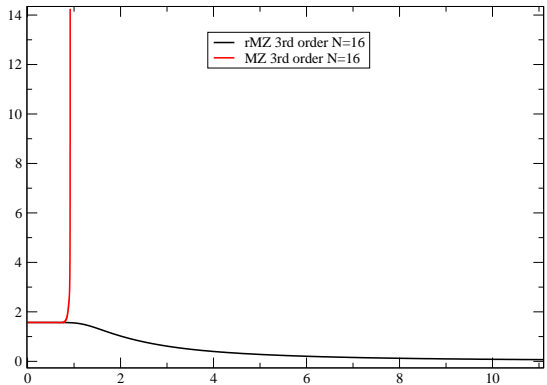,width=3.5in}
\caption{Comparison, for 1D Burgers,  of the value of the energy content of the resolved modes for the third order renormalized and unrenormalized MZ models.}
\label{plot_burgers_renovsnoreno}
\end{figure}

\begin{figure}
\centering
\epsfig{file=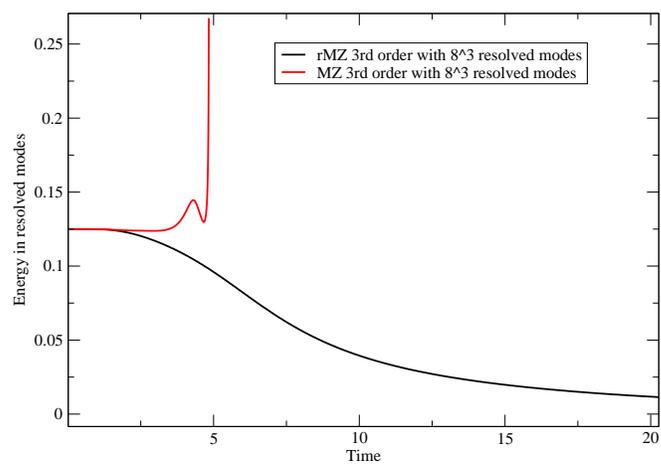,width=4.in}
\caption{Comparison, for 3D Euler,  of the value of the energy content of the resolved modes for the third order renormalized and unrenormalized MZ models.}
\label{plot_euler_renovsnoreno}
\end{figure}

In this section we present results which show that the renormalized version of the MZ formalism is advantageous with respect to the original MZ formalism. In particular, we show that for the same order in the Taylor expansion of the memory term, the renormalized algorithm leads to the stabilization of the reduced model. Figure \ref{plot_burgers_renovsnoreno} compares, for the 1D Burgers equation, the energy $\frac{1}{2}\sum_{k \in F} |u_k|^2$ for 16 resolved modes for the renormalized and unrenormalized third order models. Figure \ref{plot_euler_renovsnoreno} compares, for the 3D Euler equations, the energy $\frac{1}{2}\sum_{k \in F} |u_k|^2$ for $8^3$ resolved modes for the renormalized and unrenormalized third order models. The unrenormalized models quickly become unstable and lose all predictive ability.

The unrenormalized expansion leads to divergence of the predicted energy content of the resolved modes. This is analogous to the divergences that plagued perturbative calculations in quantum field theory (QFT) before the advent of renormalization \cite{delamotte}. In QFT, the reason for the divergences was that the perturbation expansion was performed in a quantity (the bare mass, charge etc.) which turns out to be ill-defined. The process of renormalization replaces the perturbation expansion in powers of the ill-defined quantity with a perturbation expansion in powers of the experimentally determined values of this quantity. This allows the subtraction of the terms that cause divergences and leads to finite results. 

In our case, the expansion of the memory of the MZ formalism is ill-defined because the Taylor expansion at $t=0$ breaks down after some time. On the other hand, the renormalized MZ formalism takes into account dynamic information from the evolution of the full system (while this system is still valid) and prescribes to each term in the memory expansion an appropriate coefficient. The coefficient measures how important this term is. In this way, the divergences are averted and the results become finite. 

%%%%%%%%%%%%%%End of Section Numerics%%%%%%%%%%%%%%%%%

\section{Conclusions and future work}\label{conclusions}

We have presented a new way of computing reduced models for systems of ordinary differential equations. The approach combines renormalization and effective field theory techniques with the Mori-Zwanzig formalism. The constructed reduced models are stable because they transfer activity out of the resolved scales at a rate which is dictated by the full system. The consistency between the rate of transfer activity of the reduced model and the rate of transfer activity dictated by the full system is the analog of the matching conditions employed in effective field theory. The matching conditions lead to a redefinition (renormalization) of the coefficients of a reduced model originally constructed through the Mori-Zwanzig formalism. 

The results we have obtained for the 1D Burgers and 3D Euler equations are rather encouraging. The proposed approach of calculating the renormalized coefficients leads to a linear algebraic system for the coefficients. The matrix of this linear system is ill-conditioned by construction. The reason is that the reduced model comprises of two parts. The Markovian part which is of the same functional form as the full system but defined only on the resolved scales, and the memory part which has no analog in the full system. The memory part is there to compensate for whatever activity the Markovian term cannot reproduce. Thus, the determination of the coefficients appearing in the memory part are, in essence, computed by subtracting the contribution of the Markovian term from the full system. This difference may be small and also the different terms in the memory part can have varying contributions to the difference. These factors result in the ill-conditioning of the linear system that needs to be solved. In the future, we plan to address the problem through various techniques designed to deal with ill-conditioned matrices. 

The results for 1D Burgers and 3D Euler equations also point to the possibility that the renormalized coefficients may be completely determined by the structure of the initial condition and the scaling symmetries of the equation. This needs to be investigated further and a detailed study involving random initial conditions with various spectra will be presented elsewhere.  

The proposed approach can also be applied to the Navier-Stokes equations \cite{doering}. It is straightforward to incorporate viscosity in the new expansion of the memory term proposed in the current work. We note that the viscosity starts contributing from the second order memory term. Also, the inclusion of viscosity does not complicate considerably the recursive algorithm for the calculation of the higher order terms. The expressions needed to compute the viscosity contributions can be estimated through terms already computed in the Pascal triangle construction for the inviscid terms.

The approach presented in the current work opens new possibilities for the construction of accurate and stable reduced models for (large) systems of ordinary differential equations. It also highlights the affinity between problems of model reduction in scientific computing and the construction of effective field theories in high energy physics. We hope that this connection will benefit the problem of constructing reduced models and will be of use in tackling real world problems which are impossible to address through brute force calculations. 

In conclusion, as Steven Weinberg once put it \cite{weinberg}, renormalization is indeed a good thing.

\section*{Acknowledgements} I am grateful to Profs. G.I. Barenblatt, A.J. Chorin and O.H. Hald for their ongoing guidance and support.


\begin{thebibliography}{99}

\bibitem{doering}
Doering C.R. and Gibbon J.D., Applied Analysis of the Navier-Stokes Equations, Cambridge University Press, 1995.

%\bibitem{C94}
%Chorin, A.J.,  Vorticity and Turbulence, Springer, NY, 1994.

%\bibitem{lax}
%Lax, P.D., Hyperbolic Systems of Conservation Laws and the Mathematical Theory of
%Shock Waves, SIAM Publications, Philadelphia, 1972.
   
\bibitem{sulem}
Sulem C. and Sulem P.-L., The nonlinear Schr\"odinger Equation - Self-focusing and wave collapse, Applied Mathematical Sciences, 139, Springer, New York, 1999.

\bibitem{givon}
Givon, D., Kupferman, R. and Stuart, A., Extracting macroscopic 
dynamics: model problems and algorithms, Nonlinearity 17 (2004) pp. R55-R127.

\bibitem{CS05}
Chorin, A.J. and Stinis, P., Problem reduction, renormalization and memory,
Comm. App. Math. Comp. Sci. 1 (2005) pp. 1-27.

\bibitem{CHK00}
Chorin, A.J., Hald, O.H. and Kupferman, R.,
Optimal prediction and the Mori-Zwanzig representation of irreversible
processes. Proc. Nat. Acad. Sci. USA 97 (2000) pp. 2968-2973.

\bibitem{CHK3}
Chorin, A.J., Hald, O.H. and Kupferman, R., Optimal prediction with memory, Physica D 166 (2002) pp. 239-257.

\bibitem{bernstein}
Bernstein D., Optimal prediction of Burger's equation, Multi. Mod. Sim. 6 (2007) pp. 27-52.

\bibitem{HS06}
Hald O.H. and Stinis P., Optimal prediction and the rate of decay for solutions of the Euler equations in two and three dimensions, Proc. Natl. Acad. Sci.,104, no. 16 (2007) pp. 6527-6532.

\bibitem{S10}
Stinis P., Numerical computation of solutions of the critical nonlinear Schr\"odinger equation after the singularity, arXiv:1010.2246v1 (2010).

\bibitem{chandy}
Chandy A.J. and Frankel S.H., The $t$-model as a large eddy simulation model for Navier-Stokes equations, Multi. Mod. Sim., 8, (2009) pp. 445-462.

\bibitem{georgi}
Georgi H., Effective Field Theory, Annu. Rev. Nucl. Part. Sci., Vol. 43 (1993) pp. 209-252. 

\bibitem{collins}
Collins J., Renormalization, Cambridge University Press, Cambridge, 1984.

\bibitem{goldenfeld}
Goldenfeld, N., Lectures on Phase Transitions and the Renormalization Group, Perseus Books, Reading, Mass., 1992.

\bibitem{S09}
Stinis P., A phase transition approach to detecting singularities of PDEs, Comm. Appl. Math. Comp. Sci., Vol. 4 (2009), No. 1, 217-239.

\bibitem{binney}
Binney J., Dowrick N., Fisher A., Newman M., The Theory of Critical Phenomena (An Introduction to the Renormalization Group), The Clarendon Press, Oxford, 1992.

\bibitem{golub}
Golub G.H. and Van Loan C.F., Matrix Computations, Third Edition, The Johns Hopkins University Press, Baltimore, 1996.


\bibitem{S06}
Stinis P., Higher order Mori-Zwanzig models for the Euler equations, 
Multi. Mod. Sim. 6, no 3, (2007) pp. 741-760.

\bibitem{bellman}
Bellman R., Perturbation Techniques in Mathematics, Engineering and Physics, Dover, NY, 2003. 


\bibitem{boyd}
Boyd J.P., Chebyshev and Fourier Spectral Methods, Dover, New York, 2001.

\bibitem{chorinran}
Chorin A.J., Random choice solution of hyperbolic systems, J. Comp. Phys. 22 (1976) pp. 517-533.

\bibitem{hairer}
Hairer E.,  N\"orsett S.E.  and Wanner G., Solving Ordinary Differential Equations I, Springer, NY, 1987.

\bibitem{gottlieb}
Don W.S., Gottlieb D., Shu C.W., Schilling O. and Jameson L., J. Sci. Comp. 24 (2005), pp. 569-595.

\bibitem{delamotte}
Delamotte B, A hint of renormalization, Am. J. Phys. Vol. 72, no. 2 (2004) pp. 170.

%\bibitem{majda}
%Majda, A.J. and Bertozzi, A.L., Vorticity and Incompressible Flow, Cambridge University Press, Cambridge, 2002.


%\bibitem{wilson}
%Wilson K.,  The renormalization group and critical phenomena, Rev. Mod. Phys. 55 (1983) pp. 583-600.

\bibitem{weinberg}
Weinberg S., Why the renormalization group is a good thing, Asymptotic Realms of Physics: Essays in Honor of Francis E. Low, MIT Press, Cambridge MA, 1983.





\end{thebibliography}
\end{document}